\documentclass[preprint,12pt]{elsarticle}

\usepackage{amssymb}
\usepackage{amsmath,amssymb,amsfonts}
\usepackage{tikz}

\usetikzlibrary{decorations.pathmorphing, patterns, decorations.markings}

\newtheorem{thm}{Theorem}
\newtheorem{Cor}{Corollary}
\newtheorem{Lem}{Lemma}
\newdefinition{rmk}{Remark}
\newproof{proof}{Proof}
\journal{Applied Mathematical Modelling}

\begin{document}

\begin{frontmatter}

\title{An accurate approach to determining the spatiotemporal vehicle load on bridges based on
measured boundary slopes}

\author[1]{Alemdar Hasanov\fnref{fn1}}
\ead{alemdar.hasanoglu@gmail.com}
\author[2]{Onur Baysal\corref{cor1}
\fnref{fn2}}
\ead{onur.baysal@um.edu.mt}

\cortext[cor1]{Corresponding author}
\fntext[fn1]{Emeritus Professor. \c{S}ehit Ekrem Dsitrict, Altun\c{s}ehir Str., Ayazma Villalari, No: 22. Bah\c{c}ecik - Ba\c{s}iskele, Kocaeli, 41030}
\fntext[fn2]{Senior Lecturer, Department of Mathematics, University of Malta, Malta}

\address[1]{Kocaeli University, 41001, \.{I}zmit/Kocaeli, T\"{u}rkiye}
\address[2]{University of Malta, Msida, Malta}

\begin{abstract}
In this paper, a novel mathematical model is developed to evaluate the spatiotemporal vehicle loads on long bridges from slope measurements made at the ends of a bridge based on  Euler-Bernoulli beam model with internal and external damping. The mathematical modelling of this phenomena leads to the inverse source problem of determining the spatiotemporal vehicle load $F(x,t)$ in the variable coefficient Euler-Bernoulli equation $\rho_A(x)u_{tt}+\mu(x) u_{t}+(r(x)u_{xx})_{xx}+(\kappa(x)u_{xxt})_{xx}=F(x,t)$,
$(x,t)\in \Omega_T:=(0,\ell)\times (0,T)$ subject to the "simply supported"  boundary conditions $u(0,t)=(r(x)u_{xx}+(\kappa(x)u_{xxt})_{x=0}=0$, $u(\ell,t)=(r(x)u_{xx}+(\kappa(x)u_{xxt})_{x=\ell}=0$, from the both measured outputs: $\theta_1(t):=u_x(0,t)$ and $\theta_2(t):=u_x(\ell,t)$, that is, the measured boundary slopes. It is shown that the input-output maps $(\Phi F)(t):=u_x(0,t;F)$, $(\Psi F)(t):=u_x(\ell,t;F)$, $F \in \mathcal{F}\subset L^2(\Omega_T)$, corresponding to the inverse problem, are compact and Lipschitz continuous. Then Tikhonov functional $J(F)=\Vert \Phi F-\theta_1 \Vert_{L^2(0,T)}^2+\Vert \Psi F-\theta_2 \Vert_{L^2(0,T)}^2$ is introduced to prove the existence of a quasi-solution to the inverse problem. An explicit gradient formula for the Fr\'{e}chet derivative of the Tikhonov functional is derived. The Lipschitz continuity of the Fr\'{e}chet gradient, which guarantees the monotonicity of iterations in gradient methods, has been proven.
\end{abstract}

\end{frontmatter}


\section{Introduction}

Incorrectly assessing vehicle-induced loads can lead to fatigue cracking or even bridge collapse. Furthermore, establishing relationships between these loads using experimentally measurable data is an important basis for bridge design, safety assessment, maintenance, and reinforcement. The statistical data obtained from these relationships can help us better understand the behavior of bridges under vehicle loads. Therefore, accurate and reliable estimation of vehicle loads is very important.

The dynamic effects of moving loads on bridges were not recognized until the mid-19th century. The analysis and history of moving load problems are described by Timoshenko in his book \cite{Timoshenko:1953}. The monograph \cite{Fryba:1972} was the first to give the analysis of numerous basic moving load problems and their analytical solutions. Later, these problems  extensively were reviewed in the book \cite{Yang:2004}. An overview, analysis and history of dynamic problems caused by moving loads are given in article \cite{Ouyang:2011}. An isogeometric approach to dynamic analysis of beam constructions subjected to moving vehicles was proposed in \cite{VanDo:2017}, based on the Euler-Bernoulli model.

The primary live loads on bridges are vehicle loads, which are critical parameters for bridge health monitoring. However, traditional Weigh-in-Motion (WIM) systems \cite{Lydon:2016} require the installation of weighing devices embedded in the road surface, which requires traffic interruptions during installation. However, this process is time-consuming and costly, which prevents the method from being widely implemented. As a sequel, Bridge Weigh-in-Motion (BWIM) systems \cite{Moses:1979}, developed in the 1970s, offer relatively easy installation and somewhat lower installation costs compared to traditional WIM systems. Over the last 40 years, BWIM systems have achieved a high level of research maturity, established technical methodologies, and high measurement accuracy. However, they still require special instrumentation and have complex installation requirements. Furthermore, the specificity of BWIM systems, where each bridge requires a special system, results in high initial installation and subsequent maintenance costs. This limits the deployment of dynamic weighing systems to only a few specific bridges.

In recent years, major developments in computer vision and image processing technologies have attracted a lot of attention in the field of civil engineering. In this context, the study of using computer vision to identify vehicle loads on bridges has also been proposed by various papers. In  \cite{Ojio:2016} the authors propose a method that uses a controller to simultaneously control the cameras on the bridge deck and underneath the bridge. The camera on the bridge deck can measure the
dynamic displacement response data of the bridge data and receive the vehicle's axle information. They proved the applicability of employing computer vision to determine vehicle weight by identifying vehicle loads through the analysis of both sets of data. In \cite{Martini:2022}, a computer vision-based approach to determine the tire loads of moving cars, the load places on the bridge, and the displacement response of the bridge was proposed. A bridge influence line model was then built in accordance with these findings. This technique made use of multiple cameras cooperating. By identifying tire types and getting tire pressure information from a database, the on-bridge cameras calculated vehicle loads. Vehicle load positions and the bridge displacement reaction were picked up by the under-bridge cameras. They created a bridge influence line model by integrating these data, offering a fresh method for computer vision-based bridge structure monitoring.

The use of computer vision for non-contact, target-free displacement measurements was discussed in detail in \cite{Khuc:2018}. Here, an iterative approximation algorithm to construct the displacement unit influence surface (UIS) was proposed. This algorithm allows to estimate the equivalent moving loads on the bridge under multiple vehicle loads using camera data and computer vision algorithms.

To solve for vehicle load, a spatiotemporal connection model of structural displacement, vehicle load, and load distribution was developed in \cite{Tang:2024}. To confirm the viability of the proposed method, engineering practice experiments and model bridge testing under varied loading scenarios are carried out. According to the results of the model bridge tests, the structural displacement determined by traffic video measurement can accurately represent the displacement characteristics of the structure.

As an alternative to the above methods and approaches, in this paper we propose a new mathematical model of vehicle-bridge interaction based on the damped Euler–Bernoulli beam, which takes into account all the main physical parameters, and a new solution algorithm based on weak solution theory for the forward problem and a quasi-solution approach combined with the adjoint method for the inverse problem. In this model, the angles at the ends of the bridge are used as measurable data. The advantage of this approach is that it allows for almost zero-cost load estimation on bridges equipped with traffic surveillance, meets more practical needs with extremely low monitoring costs, and promotes the widespread application of vehicle load estimation. In addition, unlike known methods, the error margin of this method in load estimation is below $10\%$.

The rest of the paper is organized as follows. In Section 2, we describe the mathematical model of bridge-vehicle load interaction, and formulate the identification problem. Analysis of the weak solution of the corresponding forward problem with a'pripori estimates is discussed in Section 3. In Section 4, we introduce the input-output operator and prove some properties of this operator.
The corresponding Tikhonov functional is introduced in Section 5. The Fr\'echet derivative of this functional through a suitable adjoint problem is derived in Section 6. In Section 7, the Lipschitz continuity of the Fr\'{e}chet gradient is proved. 


\section{The mathematical model of bridge-vehicle load interaction}

Within the Euler-Bernoulli damped beam model, the vibration of a simply supported long bridge under  the spatiotemporal load $F(x,t)$ is described by the following mathematical model:
\begin{eqnarray}\label{1}
\left\{ \begin{array}{ll}
\rho_A(x) u_{tt}+\mu(x)u_{t}-(T(x)u_{x})_{x}+ (r(x)u_{xx}+\kappa(x)u_{xxt})_{xx} =F(x,t),\\
\qquad \qquad \qquad \qquad \qquad \qquad \qquad \qquad (x,t)\in \Omega_{T}:=(0,\ell)\times (0,T);\\ [2pt]
u(x,0)=u_{t}(x,0)=0, ~x \in (0,\ell); \\ [2pt]
u(0,t)=\left(r(x)u_{xx}+\kappa(x)u_{xxt}\right)_{x=0}=0; \\
\qquad \qquad \qquad \qquad u(\ell,t)=\left (r(x)u_{xx}+\kappa(x)u_{xxt}\right)_{x=\ell}=0,~t \in [0,T].
\end{array} \right.
\end{eqnarray}

Here and below, $u(x,t)$ is the transverse deflection in position $x\in (0,\ell)$ and at the time $t\in \times [0,T]$, while $T>0$ is the final time instance, which may be small enough, and $\ell>0$ is the length of the beam. Further, $\rho_A(x):=\rho(x)A_s(x)$, while $\rho(x)>0$ and $A_s(x)>0$ are the mass density and the cross-sectional area, $r(x):=E(x)I(x)>0$ is the flexural rigidity (or bending stiffness) of a nonhomogeneous beam while $E(x)>0$ is the elasticity modulus and $I(x)>0$ is the moment of inertia. $T_r(x)\ge 0$ is the axial tensile force.

The external and internal damping mechanisms are given by the terms $\mu(x)u_t$ and $(\kappa(x)u_{xxt})_{xx}$, respectively. The coefficients $\mu(x)\ge 0$ and $\kappa(x)>0$ are called the viscous (internal) damping and the strain-rate or Kelvin-Voigt damping coefficients, respectively. The coefficient $\kappa(x):=c_d(x)I(x)$ represents energy dissipated by friction internal to the beam, while $c_d>0$ is the strain-rate damping coefficient.

The function $F(x,t)$ is the spatiotemporal load expressing the effect of the moving car on the bridge.

Figure 1 represents the geometry of the problem (\ref{1}).

Suppose that the slopes $\theta_0(t)$ and $\theta_\ell(t)$ at the ends $x=0$  and  $x=\ell$ of the bridge are given as measurable data:
\begin{eqnarray}\label{2}
\theta_0(t):=u_x(0,t),~\theta_\ell(t):=u_x(\ell,t),~ t \in [0,T].
\end{eqnarray}

Within the framework of the mathematical model (\ref{1})-(\ref{2}), the problem of determining the unknown spatiotemporal load $F(x,t)$ based on these data, is defined as follows:\\
\emph{Find the unknown spatiotemporal load $F(x,t)$ in (\ref{1}) from the measured slopes $\theta_0(t)$ and $\theta_\ell(t)$ introduced in (\ref{2}).}

The problem (\ref{1})-(\ref{2}) is defined as a spatiotemporal load identification problem or an inverse source problem with two Dirichlet measured outputs $\theta_0(t)$ and $\theta_\ell(t)$, according to generally accepted terminology.

For a given function $F(x,t)$ from some class of admissible loads, the initial boundary value problem (\ref{1}) will be referred as the \emph{forward problem}.

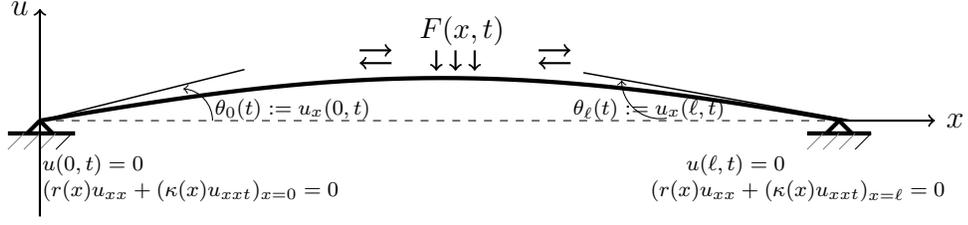
\begin{figure}
  \centering
 \hspace*{0.6cm} \begin{tikzpicture}[scale=.85]
      \draw[color=black,ultra thick] (0.0,0.0) -- (-0.2,-0.2);
      \draw[color=black,ultra thick] (0.0,0.0) -- (0.2,-0.2);
      \draw[color=black,ultra thick] (-0.50,-0.20) -- (0.55,-0.20);
      \draw[color=gray] (-0.50,-0.45) -- (-0.25,-0.20);
      \draw[color=gray] (-0.25,-0.45) -- (0,-0.20);
      \draw[color=gray] (0,-0.45) -- (0.25,-0.20);
      \draw[color=black] (0.25,-0.45) -- (0.50,-0.20);
     \draw [line width=0.2mm, black] (0.0,0.0) -- (3.2,0.8);
     \draw [line width=0.2mm, black] (8.5,0.75) -- (12.7,0.0);
      \node[label=right:{\scriptsize $ \theta_0(t):=u_x(0,t)$}] at (2.4,0.18) {};
      \node[label=right:{\scriptsize $ \theta_\ell(t):=u_x(\ell,t)$}] at (8.0,0.18) {};
    \draw[thin, ->] (2.7,0) arc (0:85:0.5);
     \draw[thin, <-] (9.1,0.63) arc (0:100:-0.6);
     \draw[dashed] (0,0) -- (4*3.14,0);
        \draw[thick, ->] (5.0,1.1)-- (5.5,1.1);
       \draw[thick, <-] (5.0,0.9)-- (5.5,0.9);
       \draw[thick, ->] (7.8,1.1)-- (8.3,1.1);
       \draw[thick, <-] (7.8,0.9)-- (8.3,0.9);
        \node[label=right:{\small$ F(x,t)$}] at (5.6,1.4) {};
         \draw[thick, ->] (6.2,1.1)-- (6.2,0.78);
         \draw[thick, ->] (6.5,1.1)-- (6.5,0.78);
         \draw[thick, ->] (6.8,1.1)-- (6.8,0.78);
      \draw[thick,->] (4*3.14,0) -- (14,0) node [right]{$x$};
      \draw[thick,->] (0,-1.5) -- (0,1.75) node [left]{$u$};
      \draw[color=black,ultra thick,smooth,domain=0:{4*3.14}] plot (\x,{sin(0.25*deg(\x))/1.5});
      \node[label=right:{\scriptsize$ u(0,t)=0$}] at (-.3,-0.7) {};
      \node[label=right:{\scriptsize$(r(x)u_{xx}+(\kappa(x)u_{xxt})_{x=0}=0$}] at (-.3,-1.1) {};
      \node[label=left:{\scriptsize$u(\ell,t)=0$}] at (12,-0.7) {};
      \node[label=left:{\scriptsize $(r(x)u_{xx}+(\kappa(x)u_{xxt})_{x=\ell}=0$}] at (14.5,-1.1) {};
     \draw[color=black,ultra thick] (12.5,0.0) -- (-0.2+12.5,-0.2);
      \draw[color=black,ultra thick] (12.5,0.0) -- (0.2+12.5,-0.2);
      \draw[color=black,ultra thick] (12.0,-0.20) -- (13,-0.20);
      \draw[color=gray] (12.5-0.50,-0.45) -- (12.5-0.25,-0.20);
      \draw[color=gray] (12.5-0.25,-0.45) -- (12.5,-0.20);
      \draw[color=gray] (12.5,-0.45) -- (12.5+0.25,-0.20);
      \draw[color=black] (12.5+0.25,-0.45) -- (12.5+0.50,-0.20);
  \end{tikzpicture}
  \caption{Geometry of the spatiotemporal vehicle load on a long bridge based on the damped Euler-Bernoulli beam model}
\label{Fig-1}
  \end{figure}

\section{Analysis of the forward problem and a'pripori estimates}

We assume that the inputs and outputs in (\ref{1}) and (\ref{2}) satisfy the following physically justified basic conditions with a minimum requirement of smoothness:
\begin{eqnarray} \label{3}
\left \{ \begin{array}{ll}
\rho_A, \mu, T_r, r, \kappa \in L^{\infty}(0,\ell),~F \in L^2(\Omega_T),\\ [1pt]
\theta_0,\theta_\ell \in L^2(0,T),\\ [1pt]
0 < \rho_0 \le \rho_A (x)\le \rho_1,~0\le T_{r0} \le T_r(x) \le T_{r1},~0< r_0 \le r(x) \le r_1, \\ [1pt]
0 \le \mu_0\le \mu(x)\le  \mu_1,~0 < \kappa_0 \le \kappa(x)\le \kappa_1,~x \in (0,\ell).
\end{array} \right.
\end{eqnarray}

We use the weak solution theory developed in \cite{Baysal:2019} and \cite{Hasanov:Romanov:2021} to derive some a'priori estimates for the weak solution Euler-Bernoulli beam equation subject to clamped boundary conditions. For cantilever beams with Kelvin-Voigt damping, similar estimates have been derived in \cite{Sakthivel:2024} .

\begin{thm}\label{Theorem-1}
Assume that the basic conditions (\ref{3}) hold. Then for the weak solution
$u\in L^2(0,T;\mathcal{V}^2(0,\ell))$ with $u_t\in L^2(0,T;L^2(0,\ell))$, $u_{tt}\in L^2(0,T;H^{-2}(0,\ell))$ of the forward problem (\ref{1}), the following estimates are satisfied:
\begin{eqnarray}\label{4}
\left. \begin{array}{ll}
\displaystyle \Vert u_t \Vert^2_{L^{\infty}(0,T;L^2(0,\ell))}
\leq \frac{1}{\rho_0}\,C_e^2\, \Vert F \Vert^2_{L^2(\Omega_T)},  \\ [10pt]
\Vert u_t \Vert^2_{L^2(0,T;L^2(0,\ell))} \leq (C_e^2 -1) \,\Vert F \Vert^2_{L^2(\Omega_T)},\\ [12pt]
\displaystyle \Vert u_{xx} \Vert^2_{L^{\infty}(0,T;L^2(0,\ell))}
\leq \frac{1}{r_0}\,C_e^2\, \Vert F \Vert^2_{L^2(\Omega_T)},  \\ [10pt]
\displaystyle \Vert u_{xx} \Vert^2_{L^2(0,T;L^2(0,\ell))}
\leq  \frac{\rho_0}{r_0}\,(C_e^2 -1)\, \Vert F \Vert^2_{L^2(\Omega_T)},  \\ [13pt]
\displaystyle \Vert u_{xxt} \Vert^2_{L^{\infty}(0,T;L^2(0,\ell))}
\leq \frac{1}{\kappa_0}\,C_e^2\, \Vert F \Vert^2_{L^2(\Omega_T)},  \\ [10pt]
\displaystyle \Vert u_{xxt} \Vert^2_{L^2(0,T;L^2(0,\ell))}
\leq \frac{\rho_0}{\kappa_0}\,(C_e^2 -1)\, \Vert F \Vert^2_{L^2(\Omega_T)},
\end{array} \right.
\end{eqnarray}
where $\mathcal{V}^2(0,\ell):=\{v\in H^2(0,\ell):~ v(0)=v(\ell)=0\}$,
\begin{eqnarray}\label{5}
C_e^2=\exp(T/\rho_0),
\end{eqnarray}
and $\rho_0,r_0>0$ are the constants introduced in (\ref{3}).
\end{thm}
{\bf Proof.} Utilize the following identities
\begin{eqnarray}\label{6}
\left. \begin{array}{ll}
\displaystyle 2\int_0^t \int_0^\ell \left (r(x)u_{xx}\right )_{xx} u_{\tau} dx d\tau \\ [6pt]
\qquad = \displaystyle 2\int_0^t \int_0^\ell \left \{\left [(r(x)u_{xx})_x u_{\tau}-r(x)u_{xx} u_{x\tau}\right ]_x+\frac{1}{2}\,\left (r(x)u_{xx}^2\right )_{\tau} \right \} dx d\tau,\\ [16pt]
\displaystyle 2\int_0^t \int_0^\ell \left (\kappa(x)u_{xx\tau}\right)_{xx} u_{\tau} dx d\tau \\ [6pt]
\qquad = \displaystyle 2\int_0^t \int_0^\ell \left \{ \left [(\kappa (x)u_{xx\tau})_x u_{\tau}
-\kappa(x)u_{xx\tau} u_{x\tau}\right ]_x+\kappa (x)u^2_{xx\tau} \right \}dx d\tau,
\end{array} \right.
\end{eqnarray}
after multiplying both sides of equation (\ref{1}) by $2u_t(x,t)$, integrating the result over $\Omega_t:=(0,\ell)\times (0,t)$, $t\in (0,T)$. In view of the homogeneous initial and boundary conditions, next we get the following \emph{energy identity}:
\begin{eqnarray}\label{7}
\int_0^\ell \left [\rho_A(x) u_t^2 +r(x)u_{xx}^2+\kappa(x)u_{xx\tau}^2 +T_r(x)u_x^2\right ]dx \qquad \qquad \qquad \quad
 \nonumber \\ [1pt]
\qquad \qquad +2 \int_0^t \int_0^\ell \mu (x)u_{\tau}^2 dx d \tau = 2\int_0^t \int_0^\ell F(x,\tau) u_{\tau}dx d \tau,~t\in[0,T].
\end{eqnarray}
The following main integral inequality is derived from the identity (\ref{7}) under the basic conditions (3):
\begin{eqnarray}\label{8}
\rho_0 \int_0^\ell u_t^2dx +r_0 \int_0^\ell u_{xx}^2dx+\kappa_0 \int_0^\ell u_{xx\tau}^2 dx
+T_{r0} \int_0^\ell u_x^2dx \qquad \qquad  \nonumber \\ [1pt]
\qquad +2 \int_0^t \int_0^\ell \mu (x) u_{\tau}^2 dx d \tau \le
 \int_0^t \int_0^\ell u^2_{\tau}dx d \tau +\int_0^t \int_0^\ell F^2(x,\tau)dx d \tau,~
\end{eqnarray}
for all $t\in [0,T]$.

The first consequence of (\ref{8}) is the following inequality:
\begin{eqnarray*}
\rho_0 \int_0^\ell u_t^2dx \le \int_0^t \int_0^\ell u^2_{\tau}dx d \tau +\int_0^t \int_0^\ell F^2(x,\tau)dx d \tau,~t\in [0,T].
\end{eqnarray*}
Here the Gr\"onwall-Bellman inequality is used to get:
\begin{eqnarray}\label{9}
\displaystyle \int_0^\ell u_t^2dx \le \frac{1}{\rho_0}\, \Vert F \Vert^2_{L^2(\Omega_T)}
\exp(t/\rho_0),~t\in [0,T].
\end{eqnarray}
The first two estimates in (\ref{4}) are derived from this inequality.

As the second consequence of the main integral inequality (\ref{8}), we deduce that
\begin{eqnarray*}
r_0 \int_0^\ell u_{xx}^2dx \le
 \int_0^t \int_0^\ell u^2_{\tau}dx d \tau +\int_0^t \int_0^\ell F^2(x,\tau)dx d \tau,~
\end{eqnarray*}
With the estimate
\begin{eqnarray}\label{10}
\displaystyle \int_0^\tau \int_0^\ell u_t^2dx \le  \Vert F \Vert^2_{L^2(\Omega_T)}
\left [\exp(t/\rho_0)-1\right ],~t\in [0,T].
\end{eqnarray}
which is obtained by integration of the inequality (\ref{9}), this inequality yields:
\begin{eqnarray*}
r_0 \int_0^\ell u_{xx}^2dx \le \Vert F \Vert^2_{L^2(\Omega_T)}\exp(t/\rho_0),~t\in [0,T].
\end{eqnarray*}
This inequality readily yields the third and fourth estimations in (\ref{4}).

Finally, we use estimate (\ref{10}) with the third consequence
\begin{eqnarray*}
\kappa_0 \int_0^\ell u_{xx\tau}^2 dx \le
 \int_0^t \int_0^\ell u^2_{\tau}dx d \tau +\int_0^t \int_0^\ell F^2(x,\tau)dx d \tau
\end{eqnarray*}
of the main integral inequality  to obtain the following estimate:
\begin{eqnarray*}
\kappa_0 \int_0^\ell u_{xx\tau}^2 dx \le
 \Vert F \Vert^2_{L^2(\Omega_T)}\exp(t/\rho_0),~t\in [0,T].
\end{eqnarray*}
The fifth and sixth estimates in (\ref{4}) are easily obtained from this inequality.
 \hfill$\Box$

\begin{rmk}\label{Remark-1}
This theorem, which obviously generalizes Theorem 11.1.2 of \cite{Hasanov:Romanov:2021} for the case when the Kelvin-Voigt damping  coefficient $\kappa(x)$ is in the equation, allows us to obtain important trace estimates. In the case when $\kappa(x)$, the fifth and sixth estimates in (\ref{4})  can be obtained only for the regular weak solution of the forward problem (\ref{1}). Therefore, these considerations can also be interpreted as the Kelvin-Voigt damping coefficient increasing the regularity of the weak solution \cite{Sakthivel:2024}.
\end{rmk}

\begin{Cor}\label{Corollary-1}
If  conditions of Theorem \ref{Theorem-1} are satisfied, then the following trace estimates are true:
\begin{eqnarray}\label{11}
\left. \begin{array}{ll}
\displaystyle \Vert u_{x}(0,\,\cdot) \Vert^2_{L^2(0,T)}
\leq  \frac{C_1^2}{r_0}\, \Vert F \Vert^2_{L^2(\Omega_T)},  \\ [10pt]
\displaystyle \Vert u_{xt}(0,\,\cdot) \Vert^2_{L^2(0,T)}
\leq  \frac{C_1^2}{\kappa_0}\, \Vert F \Vert^2_{L^2(\Omega_T)},\\ [12pt]
\displaystyle \Vert u_{x}(\ell,\,\cdot) \Vert^2_{L^2(0,T)}
\leq  \frac{C_1^2}{r_0}\, \Vert F \Vert^2_{L^2(\Omega_T)},  \\ [10pt]
\displaystyle \Vert u_{xt}(\ell,\,\cdot) \Vert^2_{L^2(0,T)}
\leq  \frac{C_1^2}{\kappa_0}\,\Vert F \Vert^2_{L^2(\Omega_T)},
\end{array} \right.
\end{eqnarray}
where
\begin{eqnarray}\label{12}
\displaystyle C_1^2=\frac{5\ell \rho_0}{3}\,(C_e^2 -1),
\end{eqnarray}
and $C_e>0$ is the constant introduced in (\ref{5}).
\end{Cor}
{\bf Proof.} With the aid of Rolle's theorem, we have the following inequality:
\begin{eqnarray}\label{13}
\displaystyle \int_0^\ell v_x^2(x) dx \le  \frac{\ell^2}{2} \int_0^\ell v_{xx}^2(x) dx,
\end{eqnarray}
for a function $v \in H^2(0,\ell)$ satisfying the conditions $v(0)=v(l)=0$.

On the other hand, the identity
\begin{eqnarray*}
\displaystyle u_x(0,t)=- \frac{1}{\ell} \int_0^\ell \left (\ell-x)u_x(x,t) \right )_{x}dx
\end{eqnarray*}
for the weak solution of the forward problem (\ref{1}) implies:
\begin{eqnarray*}
\displaystyle u_x(0,t) \le \frac{2}{\ell^2} \left (\int_0^\ell u_x(x,t) dx\right )^2+
\frac{2}{\ell^2} \left (\int_0^\ell (\ell-x)u_{xx}(x,t)dx\right )^2 \\ [8pt]
\displaystyle \qquad  \qquad  \le \frac{2}{\ell} \int_0^\ell u^2_x(x,t) dx+
\frac{2\ell}{3} \int_0^\ell u^2_{xx}(x,t)dx.
\end{eqnarray*}
Hence,
\begin{eqnarray*}
\displaystyle \int_0^T u^2_x(0,t) dx  \le \frac{2}{\ell} \int_0^T \int_0^\ell u^2_x(x,t) dx dt+
\frac{2\ell}{3} \int_0^T \int_0^\ell u^2_{xx}(x,t)dx dt.
\end{eqnarray*}
In view of (\ref{13}), applied to the weak solution of the forward problem (\ref{1}), this yields:
\begin{eqnarray}\label{14}
\displaystyle \int_0^T u^2_x(0,t) dx  \le \frac{5 \ell}{3} \int_0^T \int_0^\ell u^2_{xx}(x,t)dx.
\end{eqnarray}
In the same way, we deduce that
\begin{eqnarray}\label{15}
\displaystyle \int_0^T u^2_{xt}(0,t) dx  \le \frac{5 \ell}{3} \int_0^T \int_0^\ell u^2_{xxt}(x,t)dx.
\end{eqnarray}

It is evident that the similar to (\ref{14}) and (\ref{15}) estimations are valid also for the norms $\Vert u_{x}(\ell,\,\cdot) \Vert^2_{L^2(0,T)}$ and $\Vert u_{xt}(\ell,\,\cdot) \Vert^2_{L^2(0,T)}$, respectively. 

The desired estimations (\ref{11}) are obtained if we take into account the third, fourth, fifth and sixth estimates in (\ref{4}) in inequalities (\ref{14}) and (\ref{15}) above. \hfill$\Box$

\section{The input-output operators}

We define the \emph{set of admissible loads (inputs)} as follows:
\begin{eqnarray}\label{16}
\mathcal{F}:=\{F \in L^2(\Omega_T):\, \Vert F \Vert^2_{L^2(\Omega_T)} \leq  C_F\}.
\end{eqnarray}
where $C_F>0$ independent on $F(x,t)$ constant.

Let $F \in \mathcal{F}$ be given. Denote by $u(x,t;F)$ the unique weak solution of the
 the forward problem (\ref{1}) corresponding to this input. Then $u_x(0,t;F)$ and $u_x(\ell,t;F)$ are the outputs, and we introduce the input-output operators as follows:
\begin{eqnarray}\label{17}
\left. \begin{array}{ll}
\Phi_0(F)(t):=u_x(0,t;F),~F \in \mathcal{F},\, t\in [0,T], \\ [6pt]
\Phi_0: \mathcal{F} \subset L^2(\Omega_T) \mapsto L^2(0,T); \\[10pt]
\Phi_\ell(F)(t):=u_x(\ell,t;F),~F \in \mathcal{F},\, t\in [0,T], \\ [6pt]
\Phi_\ell : \mathcal{F} \subset L^2(\Omega_T) \mapsto L^2(0,T).
\end{array} \right.
\end{eqnarray}
In view of these operators the inverse source problem (\ref{1})-(\ref{2}) can be reformulated as the following system of operator equations:
 \begin{eqnarray}\label{18}
\left \{ \begin{array}{ll}
\Phi_0(F)(t)=\theta_0(t),~\theta_0 \in  L^2(0,T), \\ [6pt]
\Phi_\ell(F)(t)=\theta_\ell(t),~\theta_\ell \in  L^2(0,T),~F \in \mathcal{F},\, t\in [0,T].
\end{array} \right.
\end{eqnarray}

The following lemma will make it possible for us to assert the ill-posedness of the considered inverse problem.

\begin{Lem}\label{Lemma-1}
Assume that conditions of Theorem \ref{Theorem-1} are satisfied. Then the input-output operators introduced in (\ref{17}) are compact.
\end{Lem}
{\bf Proof.} Let $\{F^{(m)} \}\subset \mathcal{F}$, $m=1,2,3,\, ...$, be a sequence of admissible loads. Denote by $\{\Phi_0(F^{(m)}) \}, \{\Phi_\ell(F^{(m)}) \}\subset \mathcal{F}$ the sequences of outputs: $\Phi_0(F^{(m)})(t)=u_x(0,t;F^{(m)})$, $\Phi_\ell(F^{(m)})(t)=u_x(\ell,t;F^{(m)})$.
From the estimates (\ref{11}), it follows that the norms are bounded in the norm of the Sobolev space $H^1(0,T)$, and hence compact in $L^2(0,T)$. This means that the  input-output operators transform a bounded in $L^2(0,T)$ set $\{F^{(m)} \}$, $m=1,2,3,\, ...$, to the compact in $L^2(0,T)$ sets $\{\Phi_0(F^{(m)}) \}$, and $\{\Phi_\ell(F^{(m)}) \}$. Hence, these operators are compact. \hfill$\Box$

\begin{rmk}\label{Remark-2}
Here, unlike Lemma 11.2.1 of \cite{Hasanov:Romanov:2021}, the compactness property is obtained without assuming the existence of a \emph{regular }weak solution, thanks to the presence of the Kelving-Voigh damping term in the mathematical model (\ref{1}). This feature was first found in the study \cite{Sakthivel:2024}.
\end{rmk}

From Lemma \ref{Lemma-1}, it follows, in particular, that the \emph{inverse problem (\ref{1})-(\ref{2}) is ill-posed}.

\begin{Lem}\label{Lemma-2}
Under the conditions of Theorem \ref{Theorem-1}, the input-output operators (\ref{17}) are Lipschitz continuous, that is,
\begin{eqnarray}\label{19}
\left. \begin{array}{ll}
\displaystyle \Vert \Phi_0(F_1)-\Phi_0(F_2) \Vert_{L^2(0,T)}
\leq  C_L \Vert F_1-F_2 \Vert_{L^2(\Omega_T)},  \\ [8pt]
\displaystyle \Vert \Phi_\ell(F_1)-\Phi_\ell(F_2) \Vert_{L^2(0,T)}
\leq  C_L \Vert F_1-F_2 \Vert_{L^2(\Omega_T)}, ~\forall F_1,F_2 \in \mathcal{F},
\end{array} \right.
\end{eqnarray}
where $C_L=C_1/\sqrt{r_0} >0$ is the Lipschitz constant, and $C_1>0$ is the constant introduced in (\ref{12}).
\end{Lem}
{\bf Proof.}  Denote by $\delta u(x,t):=u(x,t;F_1) -u(x,t;F_2)$ the weak solution of the forward problem (\ref{1}) with $F(x,t)$ replaced by $\delta F(x,t):=F_1(x,t)-F_2(x,t)$, $F_1,F_2 \in \mathcal{F}$. Then, from estimates (\ref{11}) it follows that
\begin{eqnarray*}
\left. \begin{array}{ll}
\displaystyle \Vert \delta u_{x}(0,\,\cdot) \Vert^2_{L^2(0,T)}
\leq  \frac{C_1^2}{r_0}\, \Vert \delta F \Vert^2_{L^2(\Omega_T)},  \\ [10pt]
\displaystyle \Vert \delta u_{x}(\ell,\,\cdot) \Vert^2_{L^2(0,T)}
\leq  \frac{C_1^2}{r_0}\, \Vert \delta F \Vert^2_{L^2(\Omega_T)}.
\end{array} \right.
\end{eqnarray*}
By definition of the input-output operators, $\delta u_{x}(0,t)=\Phi_0(F_1)-\Phi_0(F_2)$, $\delta u_{x}(\ell,t)=\Phi_\ell (F_1)-\Phi_\ell (F_2)$, and
these estimates lead to the required statements (\ref{19}).  \hfill$\Box$

We will show that the property of the input-output operators yields the Lipschitz continuity of the Tikhonov functional, which will be introduced below.

\section{The Tikhonov functional and existence of a quasi-solution}

Due to random noise in the measured outputs $\theta_0(t)$ and $\theta_\ell(t)$, it is evident that exact equalities in the system of equations (\ref{18}) are not achievable in practice. As a consequence, one needs to introduce the Tikhonov functional
\begin{eqnarray}\label{20}
\displaystyle J(F):=\frac{1}{2} \int_0^T \left [\Phi_0(F)(t) -\theta_0(t) \right ]^2dt
+ \frac{1}{2} \int_0^T \left [\Phi_\ell(F)(t) -\theta_\ell(t) \right ]^2dt \quad \nonumber \\ [6pt]
\qquad \equiv \frac{1}{2} \int_0^T \left [u_x(0,t;F) -\theta_0(t) \right ]^2dt
+ \frac{1}{2} \int_0^T \left [u_x(\ell,t;F) -\theta_\ell(t) \right ]^2dt,
\end{eqnarray}
$F \in \mathcal{F}$, and consider the following minimization problem for this functional:
 \begin{eqnarray}\label{21}
J(F_*)=\inf_{F \in \mathcal{F}} J(F).
\end{eqnarray}

A solution of the minimization problem (\ref{21}) is called a quasi-solution of the inverse source problem (\ref{1})-(\ref{2}).

\begin{Lem}\label{Lemma-3}
Let conditions of Theorem \ref{Theorem-1} hold. Then the Tikhonov functional (\ref{20}) is  Lipschitz continuous:
\begin{eqnarray}\label{22}
\displaystyle \vert J(F_1)-J(F_2) \vert_{L^2(0,T)}
\leq  C_J \Vert F_1-F_2 \Vert_{L^2(\Omega_T)}, ~\forall F_1,F_2 \in \mathcal{F},
\end{eqnarray}
where
\begin{eqnarray*}
\displaystyle C_J=\left [\frac{2C_1 C_F}{\sqrt{r_0}} +\Vert \theta_0 \Vert_{L^2(0,T)}+\Vert \theta_\ell \Vert_{L^2(0,T)}\right ] C_L\,,
\end{eqnarray*}
is the Lipschitz constant, $C_1>0$ and $C_F>0$ are the constants introduced in (\ref{12}) and (\ref{16}), respectively.
\end{Lem}
{\bf Proof.} We can easily prove that
\begin{eqnarray*}
\displaystyle \vert J(F_1)-J(F_2) \vert_{L^2(0,T)}\quad \qquad \qquad \qquad \qquad \qquad \qquad \qquad \qquad \qquad \qquad \qquad  \\ [8pt]
\qquad  \leq \frac{1}{2} \left [ \Vert \Phi_0(F_1) \Vert_{L^2(0,T)}
+\Vert \Phi_0(F_2) \Vert_{L^2(0,T)} +2\Vert \theta_0 \Vert_{L^2(0,T)}\right ]\Vert \delta \Phi_0(F)\Vert_{L^2(0,T)} \\ [8pt]
 \qquad +\displaystyle \frac{1}{2} \left [ \Vert \Phi_\ell(F_1) \Vert_{L^2(0,T)}
+\Vert \Phi_\ell(F_2) \Vert_{L^2(0,T)} +2\Vert \theta_\ell \Vert_{L^2(0,T)}\right ] \Vert \delta \Phi_\ell(F)\Vert_{L^2(0,T)},
\end{eqnarray*}
for all $F_1,F_2 \in \mathcal{F}$, where $\Vert \delta \Phi_0(F)\Vert_{L^2(0,T)}= \Vert  \Phi_0(F_1)-\Phi_0(F_2) \Vert_{L^2(0,T)}$ and $\Vert \delta \Phi_\ell(F)\Vert_{L^2(0,T)}= \Vert  \Phi_\ell(F_1)-\Phi_\ell(F_2) \Vert_{L^2(0,T)}$, respectively. Using here (\ref{11}) and (\ref{16}) we obtain:
\begin{eqnarray*}
\displaystyle \vert J(F_1)-J(F_2) \vert_{L^2(0,T)}\quad \qquad \qquad \qquad \qquad \qquad \qquad \qquad  \\ [8pt]
\qquad  \leq \left [\frac{C_1 C_F}{\sqrt{r_0}}+ \Vert \theta_0\Vert_{L^2(0,T)}\right ]\Vert \delta \Phi_0(F)\Vert_{L^2(0,T)} \\ [8pt]
 \qquad +\displaystyle  \left [\frac{C_1 C_F}{\sqrt{r_0}}+
 \Vert \theta_\ell \Vert_{L^2(0,T)}\right ] \Vert \delta \Phi_\ell(F)\Vert_{L^2(0,T)},
\end{eqnarray*}

With the continuity property of the input-output operators, this leads to the required result  (\ref{22}). \hfill$\Box$

\begin{thm}\label{Theorem-2}
Assume that the basic conditions (\ref{3}) hold, and $F \in \mathcal{F}$, where $\mathcal{F}$ is the set of admissible loads. Then the minimization problem (\ref{21}) for  the Tikhonov functional has at least one solution in $\mathcal{F}$, that is the inverse source problem -(\ref{1})-(\ref{2}) has a quasi-solution.
\end{thm}

The proof of this theorem follows the same procedure as that of Theorem 10.1.11 of (\ref{21}) \cite{Hasanov:Romanov:2021}.

\section{Fr\'echet derivative of the Tikhonov functional}

Let $F, F+\delta F \in \mathcal{F}$. Denote  by $\delta J(F):= J(F+\delta F )-J(F)$ the increment of the Tikhonov functional introduced in (\ref{20}). Then,
\begin{eqnarray}\label{23}
\displaystyle \delta J(F)= \int_0^T \left [u_x(0,t;F) -\theta_0(t) \right ]\delta u_x(0,t;F) dt +\frac{1}{2} \int_0^T \left (\delta u_x(0,t;F) \right )^2 dt \nonumber \\ [8pt]
\quad \displaystyle + \int_0^T \left [u_x(\ell,t;F) -\theta_\ell(t) \right ] \delta u_x(\ell,t;F) dt +\frac{1}{2} \int_0^T \left (\delta u_x(\ell,t;F) \right )^2 dt,~~
\end{eqnarray}

\begin{Lem}\label{Lemma-4}
Let conditions the basic conditions (\ref{3}) hold. Then the following integral relationship
holds:
\begin{eqnarray}\label{24}
\displaystyle \int_0^T p(t)\delta u_x(0,t)dt +\int_0^T q(t)\delta u_x(\ell,t)dt
= \int_0^T \int_0^\ell  \delta F(x,t) \varphi (x,t) dx dt,
\end{eqnarray}
for all $p,q \in L^2(0,T)$, where the function $\varphi (x,t)$ is the weak solution of the following backward problem:
\begin{eqnarray}\label{25}
\left\{ \begin{array}{ll}
\rho_A(x) \varphi_{tt}-\mu(x)\varphi_{t}-(T(x)\varphi_{x})_{x}+ (r(x)\varphi_{xx}-\kappa(x)\varphi_{xxt})_{xx} =0,\\ [2pt]
\qquad \qquad \qquad \qquad \qquad \qquad \qquad (x,t)\in \Omega_{T}:=(0,\ell)\times (0,T);\\ [4pt]
\varphi(x,T)=\varphi_{t}(x,T)=0, ~x \in (0,\ell); \\ [4pt]
\varphi(0,t)=0,~-\left(r(x)\varphi_{xx}-\kappa(x)\varphi_{xxt}\right)_{x=0}=p(t); \\ [2pt]
\qquad \qquad \varphi(\ell,t)=0,~\left (r(x)\varphi_{xx}-\kappa(x)\varphi_{xxt}\right)_{x=\ell}=q(t),~t \in [0,T],
\end{array} \right.
\end{eqnarray}
with the inputs $p,q \in L^2(0,T)$, and $\delta u(x,t)$  is the weak solution of the forward problem (\ref{1}) with $F(x,t)$ replaced by $\delta F(x,t)$.
\end{Lem}
{\bf Proof.} Apply the integration by parts formula several times after multiplying both sides of equation (\ref{1}) for $ \delta u (x,t)$ by an arbitrary function $\varphi(x,t)$ and integrating over  $(0,T)$. Next, we get:
\begin{eqnarray*}
\int_0^T \int_0^{\ell} \left [\rho_A(x)\varphi_{tt}-\mu(x)\varphi_{t}+\left (r(x)\varphi_{xx}- \kappa (x)\varphi_{xxt} \right)_{xx}- \left (T_r(x)\varphi_{x}\right)_{x}\right] \delta u \,dx dt ~\nonumber \\ [2pt]
\quad + \int_0^\ell \left [\rho_A(x) \delta u_t \varphi-\rho_A(x) \delta u \varphi_t +\mu(x) \delta u \varphi - \delta u \left(\kappa(x)\varphi_{xx}\right)_{xx}\right ]_{t=0}^{t=T} \,dx \nonumber \\ [2pt]
~~ + \int_0^T \left [\left(r(x) \delta u_{xx}\right)_{x}\varphi-r(x)\delta u_{xx} \varphi_{x}+
r(x)\delta u_{x} \varphi_{xx}-\delta u \left (r(x)\varphi_{xx}\right)_{x} \right ]_{x=0}^{x=\ell} \,dt \nonumber \\ [2pt]
~+ \int_0^T \left [\left(\kappa(x) \delta u_{xxt}\right)_{x}\varphi-\kappa(x)\delta u_{xxt} \varphi_{x}-
\kappa(x)\delta u_{x} \varphi_{xxt}+\delta u \left (\kappa(x)\varphi_{xxt}\right)_{x} \right ]_{x=0}^{x=\ell} \,dt \nonumber \\ [2pt]
\quad + \int_0^T \left [T_r(x) \varphi_x \delta u-T_r(x)\varphi \delta u_{x} \right ]_{x=0}^{x=\ell} \,dt = \int_0^T \int_0^\ell  \delta F(x,t) \varphi (x,t) dx dt.
\end{eqnarray*}
Here, we assume that the function $\varphi(x,t)$ solves the backward problem (\ref{25}).
In view of the homogeneous boundary, initial, and final conditions in (\ref{1}) and (\ref{25}), the required integral relationship (\ref{24}) is then obtained. \hfill$\Box$

Considering the difference between the outputs $u_x(0,t;F)$ and $u_x(\ell,t;F)$, corresponding to the admissible $F \in \mathcal{F}$, and the measured outputs $\theta_0(t)$ and $\theta_\ell(t)$ in the Tikhonov functional (\ref{20}), we choose the inputs to the backward problem (\ref{25}) based on these differences, as follows:
\begin{eqnarray}\label{26}
\left\{ \begin{array}{ll}
p(t)=u_x(0,t;F)-\theta_0(t),\\ [2pt]
q(t)=u_x(\ell,t;F)-\theta_\ell(t),~t \in [0,T],\, F \in \mathcal{F}.
\end{array} \right.
\end{eqnarray}

The backward problem (\ref{25}) with the inputs (\ref{26}) is defined as the \emph{adjoint problem, corresponding to the inverse problem (\ref{1})-(\ref{2}):}
\begin{eqnarray}\label{27}
\left\{ \begin{array}{ll}
\rho_A(x) \varphi_{tt}-\mu(x)\varphi_{t}-(T(x)\varphi_{x})_{x}+ (r(x)\varphi_{xx}-\kappa(x)\varphi_{xxt})_{xx} =0,\\ [2pt]
\qquad \qquad \qquad \qquad \qquad \qquad \qquad (x,t)\in \Omega_{T}:=(0,\ell)\times (0,T);\\ [4pt]
\varphi(x,T)=\varphi_{t}(x,T)=0, \, x \in (0,\ell); \\ [4pt]
\varphi(0,t)=0,\,\left(r(x)\varphi_{xx}-\kappa(x)\varphi_{xxt}\right)_{x=0}=u_x(0,t;F)-\theta_0(t); \\ [2pt]
\quad  \varphi(\ell,t)=0,\,\left (r(x)\varphi_{xx}-\kappa(x)\varphi_{xxt}\right)_{x=\ell}=u_x(\ell,t;F)-\theta_\ell(t),\,t \in [0,T].
\end{array} \right.
\end{eqnarray}

The adjoint problem (\ref{27}), as well as the backward problem (\ref{25}), are well-posed problems, as the change of the variable $\tau=T-t$ shows.

Furthermore, the integral relationship (\ref{24}) with the inputs introduced in (\ref{26}) is defined as the \emph{input-output relationship}:
holds:
\begin{eqnarray}\label{28}
\displaystyle \int_0^T \left [u_x(0,t;F)-\theta_0(t) \right ]\delta u_x(0,t)dt
+\int_0^T \left [u_x(\ell,t;F)-\theta_\ell(t) \right ]\delta u_x(\ell,t)dt \nonumber \\ [2pt]
\qquad \qquad \qquad = \int_0^T \int_0^\ell  \delta F(x,t) \varphi (x,t;F) dx dt.
\end{eqnarray}
This integral identity expresses, as its name suggests, a relationship between the input $F(x,t)$, the outputs $u_x(0,t;F)$, $u_x(\ell,t;F)$, and the measured outputs $\theta_0(t)$, $\theta_\ell(t)$ of the inverse problem through the solution $\varphi (x,t;F)$ of the adjoint problem (\ref{27}).

As a consequence of the increment formula (\ref{23}) for the Tikhonov functional and the input-output relationship (\ref{28}), we obtain the following formula
\begin{eqnarray}\label{29}
\displaystyle \delta J(F)=  \int_0^T \int_0^\ell \varphi (x,t;F)\,\delta F(x,t)dx dt \qquad \qquad \qquad \qquad \qquad  ~ \nonumber \\ [2pt]
~\qquad \qquad +\frac{1}{2} \int_0^T \left (\delta u_x(0,t;F) \right )^2 dt +\frac{1}{2} \int_0^T \left (\delta u_x(\ell,t;F) \right )^2 dt,~~
\end{eqnarray}
which is implies \emph{the formal gradient formula}:
\begin{eqnarray}\label{30}
\displaystyle J'(F)=  \varphi (x,t;F),\, \mbox {a.e.} ~(x,t)\in \Omega_T,\, F \in \mathcal{F}.
\end{eqnarray}

To justify the gradient formula (\ref{30}), we need to examine the adjoint problem (\ref{27}) in more detail, obtaining also the necessary estimates for the weak solution. As we will see below, this formula may not be valid for every measured output data $\theta_0, \theta_\ell \in L^2(0,T)$, which are, at the same time, inputs for the adjoıint problem (\ref{27}).

We employ the change of variable $\tau=T-t$, $t\in [0,T]$, to transform the backward problem (\ref{25}) for $\varphi (x,t) \equiv \varphi (x,t;F)$ to the following initial boundary value problem
\begin{eqnarray}\label{31}
\left\{ \begin{array}{ll}
\rho_A(x) \phi_{\tau \tau}-\mu(x)\phi_{\tau}-(T(x)\phi_{x})_{x}+ (r(x)\phi_{xx}-\kappa(x)\phi_{xx\tau})_{xx} =0,\\ [2pt]
\qquad \qquad \qquad \qquad \qquad \qquad \qquad (x,\tau)\in \Omega_{T}:=(0,\ell)\times (0,T);\\ [4pt]
\phi(x,0)=\phi_{\tau}(x,0)=0, ~x \in (0,\ell); \\ [4pt]
\phi(0,\tau)=0,~-\left(r(x)\phi_{xx}-\kappa(x)\phi_{xx\tau}\right)_{x=0}=\widetilde{p}(\tau); \\ [2pt]
\qquad \qquad \varphi(\ell,\tau)=0,~\left (r(x)\phi_{xx}-\kappa(x)\phi_{xx\tau}\right)_{x=\ell}=\widetilde{q}(\tau),~\tau \in [0,T],
\end{array} \right.
\end{eqnarray}
for the function $\phi (x,\tau) = \varphi (x,t)$, with the inputs $\widetilde{p}(\tau)=p(t)$ and $\widetilde{q}(\tau)=q(t)$.

\begin{Lem}\label{Lemma-5}
Assume that in addition to the basic conditions (\ref{3}), the inputs $\widetilde{p}(\tau)$ and $\widetilde{q}(\tau)$ in (\ref{31}) satify the following regularity conditions:
\begin{eqnarray}\label{32}
\left. \begin{array}{ll}
\widetilde{p},\widetilde{q}  \in H^1(0,T).
\end{array} \right.
\end{eqnarray}
Then for the weak solution $\phi\in L^2(0,T;\mathcal{V}^2(0,\ell))$ with $\phi_\tau\in L^2(0,T;L^2(0,\ell))$, $\phi_{\tau\tau}\in L^2(0,T;H^{-2}(0,\ell))$ of the transformed problem (\ref{31}), the following estimates hold:
\begin{eqnarray}\label{33}
\left. \begin{array}{ll}
\displaystyle \Vert \phi_{xx} \Vert^2_{L^{\infty}(0,T;L^2(0,\ell))}
\leq \exp (T)\, C_0^2 \, \Vert \widetilde{Q}\,' \Vert^2_{L^2(0,T)},  \\ [8pt]
\displaystyle \Vert \phi_{xx} \Vert^2_{L^2(0,T;L^2(0,\ell))}
\leq \left (\exp (T)-1\right ) C_0^2 \, \Vert \widetilde{Q}\,' \Vert^2_{L^2(0,T)},  \\ [10pt]
\displaystyle \Vert \phi_\tau \Vert^2_{L^{\infty}(0,T;L^2(0,\ell))}
\leq \frac{\exp (T)r_0}{2\rho_0}\,C_0^2 \, \Vert \widetilde{Q}\,' \Vert^2_{L^2(0,T)},  \\ [12pt]
\displaystyle \Vert \phi_\tau \Vert^2_{L^2(0,T;L^2(0,\ell))} \leq
\frac{\left (\exp (T)-1\right )r_0}{2\rho_0}\,C_0^2 \, \Vert \widetilde{Q}\,' \Vert^2_{L^2(0,T)}, \\ [16pt]
\displaystyle \Vert \phi_{xx\tau} \Vert^2_{L^{\infty}(0,T;L^2(0,\ell))}
\leq \frac{\exp (T)r_0}{2\kappa_0}\,C_0^2 \, \Vert \widetilde{Q}\,' \Vert^2_{L^2(0,T)},  \\ [12pt]
\displaystyle \Vert \phi_{xx\tau} \Vert^2_{L^2(0,T;L^2(0,\ell))}
\leq \frac{\left (\exp (T)-1\right )r_0 }{2\kappa_0}\,C_0^2 \, \Vert \widetilde{Q}\,' \Vert^2_{L^2(0,T)},
\end{array} \right.
\end{eqnarray}
where
\begin{eqnarray}\label{34}
\left. \begin{array}{ll}
\displaystyle C_0^2= \frac{20 \ell \,C_T}{3r^2_0},~ C_T=\max (2/T,\, 1+3T/3), \\ [13pt]
\Vert \widetilde{Q} \Vert^2_{L^2(0,T)}= \Vert \widetilde{p} \Vert^2_{L^2(0,T)}+\Vert \widetilde{q} \Vert^2_{L^2(0,T)}.
\end{array} \right.
\end{eqnarray}
\end{Lem}
{\bf Proof.}
Multiply both sides of equation (\ref{31}) by $2\phi_t(x,\tau)$, integrate it over $\Omega_\tau:=(0,\ell)\times (0,\tau)$, $\tau \in [0,T]$, and use the identities (\ref{6}). Applying the integration by parts formula, using the initial and boundary conditions, after elementary transformations we obtain:
\begin{eqnarray}\label{35}
\int_0^\ell \left [\rho_A(x) \phi_\tau^2 +r(x)\phi_{xx}^2+\kappa(x)\phi_{xx\tau}^2 +T_r(x)\phi_x^2\right ]dx+2 \int_0^\tau \int_0^\ell \mu (x)\phi_{\eta}^2 dx d \eta
 \nonumber \\ [1pt]
 =2 \int_0^\tau \widetilde{p}(\eta)\phi_{x\eta}(0,\eta) d\eta+ 2 \int_0^\tau \widetilde{q}(\eta)\phi_{x\eta}(\ell,\eta) d\eta,~\tau \in[0,T].~
\end{eqnarray}
Using the $\varepsilon$-inequality and the consequence $\phi_{x}(0,0)=\phi_{x}(\ell,0)=0$ of the homogeneous initial conditions, we estimate the right-hand side integrals, as follows:
\begin{eqnarray*}
\left. \begin{array}{ll}
\displaystyle 2 \int_0^\tau \widetilde{p}(\eta)\phi_{x\eta}(0,\eta) d\eta=
-2 \int_0^\tau \widetilde{p}\,'(\eta)\phi_{x}(0,\eta) d\eta+2 \widetilde{p}(\tau)\phi_{x}(0,\tau) \qquad \qquad\\ [8pt]
\displaystyle \qquad \le \varepsilon \left [\int_0^\tau \phi^2_{x}(0,\eta) d\eta +\phi^2_{x}(0,\tau)\right ]+ \frac{1}{\varepsilon}\, \left [\int_0^T \left (\widetilde{p}\,'(\tau)\right)^2 d\tau + \left (\widetilde{p}(\tau)\right)^2\right ],~\\ [14pt]
\displaystyle 2 \int_0^\tau \widetilde{q}(\eta)\phi_{x\eta}(\ell,\eta) d\eta=-2 \int_0^\tau \widetilde{q}\,'(\eta)\phi_{x}(\ell,\eta) d\eta+2 \widetilde{q}(\tau)\phi_{x}(\ell,\tau) \qquad \qquad \\ [8pt]
\displaystyle \qquad \le \varepsilon \left [\int_0^\tau \phi^2_{x}(\ell,\eta) d\eta +\phi^2_{x}(\ell,\tau)\right ]+ \frac{1}{\varepsilon}\, \left [\int_0^T \left (\widetilde{q}\,'(\tau)\right)^2 d\tau + \left (\widetilde{q}(\tau)\right)^2\right ],
 \end{array} \right.
\end{eqnarray*}
for all $\tau \in[0,T]$. Here we use the identities
\begin{eqnarray*}
\displaystyle \widetilde{p}(\tau)=\frac {1}{\tau}\int_0^\tau \left (\eta \widetilde{p}(\eta)\right)' d\eta,~
\displaystyle \displaystyle \widetilde{q}(\tau)=\frac {1}{\tau}\int_0^\tau \left (\eta \widetilde{q}(\eta)\right)' d\eta,~\tau \in(0,T]
\end{eqnarray*}
to deduce that
\begin{eqnarray*}
\displaystyle \left (\eta \widetilde{p}(\eta)\right)^2
\le \frac {2}{\tau} \int_0^\tau \left (\widetilde{p}(\eta)\right)^2 d\eta+
\frac {2\tau}{3} \int_0^\tau \left (\widetilde{p}\,'(\eta)\right)^2 d\eta,\\ [1pt]
\displaystyle \left (\eta \widetilde{q}(\eta)\right)^2
\le \frac {2}{\tau} \int_0^\tau \left (\widetilde{q}(\eta)\right)^2 d\eta+
\frac {2\tau}{3} \int_0^\tau \left (\widetilde{q}\,'(\eta)\right)^2 d\eta,
\end{eqnarray*}
Hence,
\begin{eqnarray}\label{36}
\left. \begin{array}{ll}
\displaystyle 2 \int_0^\tau \widetilde{p}(\eta)\phi_{x\eta}(0,\eta) d\eta
\le \varepsilon \left [\int_0^\tau \phi^2_{x}(0,\eta) d\eta +\phi^2_{x}(0,\tau)\right ]\\ [8pt]
\displaystyle \qquad \qquad + \frac{1}{\varepsilon}\, \left [\frac {2}{\tau} \int_0^\tau \left (\widetilde{p}(\eta)\right)^2 d\eta+
\left (1+\frac {2\tau}{3}\right ) \int_0^\tau \left (\widetilde{p}\,'(\eta)\right)^2 d\eta\right ],~\\ [14pt]
\displaystyle 2 \int_0^\tau \widetilde{q}(\eta)\phi_{x\eta}(\ell,\eta) d\eta \le \varepsilon \left [\int_0^\tau \phi^2_{x}(\ell,\eta) d\eta +\phi^2_{x}(\ell,\tau)\right ] \\ [8pt]
\displaystyle \qquad \qquad + \frac{1}{\varepsilon}\, \left [\frac {2}{\tau} \int_0^\tau \left (\widetilde{q}(\eta)\right)^2 d\eta+
\left (1+\frac {2\tau}{3}\right ) \int_0^\tau \left (\widetilde{q}\,'(\eta)\right)^2 d\eta\right ].
 \end{array} \right.
\end{eqnarray}
For the terms in the first right-hand side square brackets, we use the trace estimate (\ref{14}) and its analogue for $x=\ell$. Then we get:
\begin{eqnarray}\label{37}
\left. \begin{array}{ll}
\displaystyle \int_0^\tau \phi^2_{x}(0,\eta) d\eta +\phi^2_{x}(0,\tau) \le \frac{5\ell}{3}\left [ \int_0^\tau \int_0^\ell \phi^2_{x}(x,\eta) d\eta dx +\int_0^\ell \phi^2_{x}(x,\tau) dx \right ],\\ [14pt]
\displaystyle \int_0^\tau \phi^2_{x}(\ell,\eta) d\eta +\phi^2_{x}(\ell,\tau) \le \frac{5\ell}{3}\left [ \int_0^\tau \int_0^\ell \phi^2_{x}(x,\eta) d\eta dx +\int_0^\ell \phi^2_{x}(x,\tau) dx\right ],
 \end{array} \right.
\end{eqnarray}
for all $\tau \in[0,T]$. Further, we estimate the terms in the second right-hand side square brackets of (\ref{36}) as follows:
\begin{eqnarray}\label{38}
\left. \begin{array}{ll}
\displaystyle \frac {2}{\tau} \int_0^\tau \left (\widetilde{p}(\eta)\right)^2 d\eta+
\left (1+\frac {2\tau}{3}\right ) \int_0^\tau \left (\widetilde{p}\,'(\eta)\right)^2 d\eta\le C_T\int_0^T \left (\widetilde{p}\,'(\eta)\right)^2 d\eta,~\\ [10pt]
\displaystyle \frac {2}{\tau} \int_0^\tau \left (\widetilde{q}(\eta)\right)^2 d\eta+
\left (1+\frac {2\tau}{3}\right ) \int_0^\tau \left (\widetilde{q}\,'(\eta)\right)^2 d\eta \le C_T \int_0^T \left (\widetilde{q}\,'(\eta)\right)^2 d\eta,
 \end{array} \right.
\end{eqnarray}
$ \tau \in[0,T]$, with $C_T>0$ introduced in (\ref{34}). Substituting (\ref{36}) with (\ref{37}) and (\ref{38}) in (\ref{35}), we deduce that 
\begin{eqnarray}\label{39}
\displaystyle \int_0^\ell \rho_A(x) \phi_\tau^2dx +\left (r_0-\frac{5\ell \varepsilon}{3}\right)\int_0^\ell \phi_{xx}^2 dx+\kappa_0 \int_0^\ell \phi_{xx\tau}^2 dx
\qquad \qquad \quad \nonumber \\ [1pt]
\qquad +\int_0^\ell T_r(x)\phi_x^2dx +2 \int_0^\tau \int_0^\ell \mu (x)\phi_{\eta}^2 dx d \eta  \qquad \qquad \qquad \qquad \qquad \nonumber \\ [1pt]
\qquad \qquad \le  \frac{5\ell\varepsilon}{3} \int_0^\tau \int_0^\ell\phi_{xx}^2 dx d \eta
 + \frac{C_T}{\varepsilon} \left [\Vert \widetilde{p}\,' \Vert^2_{L^2(0,T)}+\Vert \widetilde{q}\,' \Vert^2_{L^2(0,T)} \right ],
\end{eqnarray}
for all $\tau \in[0,T]$. We choose the arbitrary parameter $\varepsilon>0$ from the condition
$r_0-5\ell \varepsilon/3>0$, as follows:
\begin{eqnarray*}
\varepsilon= \frac{3r_0}{10\ell}.
\end{eqnarray*}
Then (\ref{39}) yields:
\begin{eqnarray}\label{40}
\displaystyle \int_0^\ell \rho_A(x) \phi_\tau^2dx +\frac{r_0}{2}\int_0^\ell \phi_{xx}^2 dx+\kappa_0 \int_0^\ell \phi_{xx\tau}^2 dx
\qquad \qquad \qquad \qquad \qquad \nonumber \\ [1pt]
 +\int_0^\ell T_r(x)\phi_x^2dx +2 \int_0^\tau \int_0^\ell \mu (x)\phi_{\eta}^2 dx d \eta  \qquad \qquad \qquad \qquad \qquad \qquad \nonumber \\ [1pt]
 \qquad \le  \frac{r_0}{2} \int_0^\tau \int_0^\ell \phi_{xx}^2 dx d \eta
 + \frac{10\ell C_T}{3r_0} \left [\Vert \widetilde{p}\,' \Vert^2_{L^2(0,T)}+\Vert \widetilde{q}\,' \Vert^2_{L^2(0,T)} \right ],
\end{eqnarray}
for all $\tau \in[0,T]$.

The required estimates (\ref{33}) are derived from the main inequality (\ref{40}), as in the proof of Theorem \ref{Theorem-1}. \hfill$\Box$

In view of (\ref{26}),
\begin{eqnarray*}
\left. \begin{array}{ll}
\Vert \widetilde{p}\,' \Vert^2_{L^2(0,T)}=\Vert p' \Vert^2_{L^2(0,T)}=
\Vert u_{xt}(0,\,\cdot) -\theta'_0\Vert^2_{L^2(0,T)},  \\ [8pt]
\Vert \widetilde{q}\,' \Vert^2_{L^2(0,T)}=\Vert q' \Vert^2_{L^2(0,T)}=
\Vert u_{xt}(\ell,\,\cdot) -\theta'_\ell \Vert^2_{L^2(0,T)}.
\end{array} \right.
\end{eqnarray*}
With the inequalities
\begin{eqnarray*}
\left. \begin{array}{ll}
\Vert u_{xt}(0,\,\cdot) -\theta'_0\Vert^2_{L^2(0,T)}\le
2\Vert u_{xt}(0,\,\cdot) \Vert^2_{L^2(0,T)}+ 2\Vert \theta'_0\Vert^2_{L^2(0,T)}, \\ [8pt]
\Vert u_{xt}(\ell,\,\cdot) -\theta'_\ell \Vert^2_{L^2(0,T)}\le
2\Vert u_{xt}(\ell,\,\cdot)\Vert^2_{L^2(0,T)}+2\Vert \theta'_\ell \Vert^2_{L^2(0,T)},
\end{array} \right.
\end{eqnarray*}
and the trace estimates in (\ref{11}) this  leads to the following estimates:
\begin{eqnarray}\label{41}
\left. \begin{array}{ll}
\displaystyle \Vert p' \Vert^2_{L^2(0,T)}\le \frac{2C_1^2}{\kappa_0} \Vert F\Vert^2_{L^2(\Omega_T)}+ 2\Vert \theta'_0\Vert^2_{L^2(0,T)}, \\ [10pt]
\displaystyle \Vert q' \Vert^2_{L^2(0,T)}\le \frac{2C_1^2}{\kappa_0} \Vert F\Vert^2_{L^2(\Omega_T)}+ 2\Vert \theta'_\ell \Vert^2_{L^2(0,T)}.
\end{array} \right.
\end{eqnarray}

The right-hand side norms $\Vert \theta'_0\Vert_{L^2(0,T)}$ and $\Vert \theta'_\ell\Vert_{L^2(0,T)}$ in estimates (\ref{41}) provide insight into the necessary conditions for existence of the weak solution to the adjoint problem (\ref{27}). Namely, the  measured outputs $\theta_0(t)$ and $\theta_\ell(t)$ must not be from the space $L^2(0,T)$, but
from the space of smoother functions $H^1(0,T)$.

\begin{thm}\label{Theorem-3}
Assume that the basic conditions (\ref{3}) are satisfied. Suppose, in addition, that the measured outputs satisfy the regularity conditions (\ref{32}). Then for there exists a weak solution $\varphi \in L^2(0,T;\mathcal{V}^2(0,\ell))$ with $\varphi_t\in L^2(0,T;L^2(0,\ell))$, $\varphi_{tt}\in L^2(0,T;H^{-2}(0,\ell))$ of the adjoint problem (\ref{27}), and the following a'piori estimates hold:
\begin{eqnarray}\label{42}
\left. \begin{array}{ll}
\displaystyle \Vert \varphi_{xx}\Vert^2_{L^{\infty}(0,T;L^2(0,\ell))}
\leq 2\exp(T) C_0^2\left [\frac{2C_1^2}{\kappa_0}\, \Vert F \Vert^2_{L^2(\Omega_T)}
+ \Vert \Theta \Vert^2_{L^2(\Omega_T)} \right ],  \\ [10pt]
\displaystyle \Vert \varphi_{xx} \Vert^2_{L^2(0,T;L^2(0,\ell))} \leq 2\left (\exp(T)-1\right) C_0^2 \left [\frac{2C_1^2}{\kappa_0}\, \Vert F \Vert^2_{L^2(\Omega_T)}
+ \Vert \Theta \Vert^2_{L^2(\Omega_T)} \right ],\\ [12pt]
\displaystyle \Vert \varphi_{t} \Vert^2_{L^{\infty}(0,T;L^2(0,\ell))}
\leq \frac{\exp(T)r_0}{\rho_0}\,C_0^2 \left [\frac{2C_1^2}{\kappa_0}\, \Vert F \Vert^2_{L^2(\Omega_T)} + \Vert \Theta \Vert^2_{L^2(\Omega_T)} \right ],  \\ [10pt]
\displaystyle \Vert \varphi_{t} \Vert^2_{L^2(0,T;L^2(0,\ell))}
\leq  \frac{\left (\exp(T)-1\right )r_0}{\rho_0}\,C_0^2 \left [\frac{2C_1^2}{\kappa_0}\, \Vert F \Vert^2_{L^2(\Omega_T)} + \Vert \Theta \Vert^2_{L^2(\Omega_T)} \right ],  \\ [13pt]
\displaystyle \Vert \varphi_{xxt} \Vert^2_{L^{\infty}(0,T;L^2(0,\ell))}
\leq \frac{\exp(T)r_0}{\kappa_0}\,C_0^2 \left [\frac{2C_1^2}{\kappa_0}\, \Vert F \Vert^2_{L^2(\Omega_T)} + \Vert \Theta \Vert^2_{L^2(\Omega_T)} \right ],  \\ [10pt]
\displaystyle \Vert \varphi_{xxt} \Vert^2_{L^2(0,T;L^2(0,\ell))}
\leq \frac{\left (\exp(T)-1\right )r_0}{\kappa_0}\,C_0^2 \left [\frac{2C_1^2}{\kappa_0}\, \Vert F \Vert^2_{L^2(\Omega_T)} + \Vert \Theta \Vert^2_{L^2(\Omega_T)} \right ],
\end{array} \right.
\end{eqnarray}
where $\mathcal{V}^2(0,\ell)$ is the subspace of the Sobolev space $H^2(0,\ell)$ introduced in Theorem \ref{Theorem-1},
\begin{eqnarray}\label{43}
\Vert \Theta \Vert^2_{L^2(\Omega_T)}= \Vert \theta_0 \Vert^2_{L^2(\Omega_T)}+ \Vert \theta_\ell \Vert^2_{L^2(\Omega_T)},
\end{eqnarray}
and $C_0,C_1>0$ are the constants defined in (\ref{34}) and (\ref{12}), respectively.
\end{thm}
{\bf Proof.} The proof of the existence of the weak solution can be done in a similar way to the proof of the related theorems in \cite{Baysal:2019} and \cite{Sakthivel:2024}.

Estimates (\ref{33}) clearly hold for the corresponding norms of the weak solution $\varphi(x,t)$ of the adjoint problem (\ref{27}), with $\Vert \widetilde {Q}\,' \Vert^2_{L^2(0,T)}$ replaced by
\begin{eqnarray*}
\displaystyle \Vert Q' \Vert^2_{L^2(0,T)}= \Vert p' \Vert^2_{L^2(0,T)}+ \Vert q' \Vert^2_{L^2(0,T)}\,.
\end{eqnarray*}
Namely,
\begin{eqnarray*}
\left. \begin{array}{ll}
\displaystyle \Vert \varphi_{xx} \Vert^2_{L^{\infty}(0,T;L^2(0,\ell))}
\leq \exp (T)\, C_0^2 \, \Vert Q' \Vert^2_{L^2(0,T)},  \\ [8pt]
\displaystyle \Vert \varphi_{xx} \Vert^2_{L^2(0,T;L^2(0,\ell))}
\leq \left (\exp (T)-1\right ) C_0^2 \, \Vert Q' \Vert^2_{L^2(0,T)},  \\ [10pt]
\displaystyle \Vert \varphi_\tau \Vert^2_{L^{\infty}(0,T;L^2(0,\ell))}
\leq \frac{\exp (T)r_0}{2\rho_0}\,C_0^2 \, \Vert Q' \Vert^2_{L^2(0,T)},  \\ [8pt]
\displaystyle \Vert \varphi_\tau \Vert^2_{L^2(0,T;L^2(0,\ell))} \leq
\frac{\left (\exp (T)-1\right )r_0}{2\rho_0}\,C_0^2 \, \Vert Q' \Vert^2_{L^2(0,T)}, \\ [10pt]
\displaystyle \Vert \varphi_{xx\tau} \Vert^2_{L^{\infty}(0,T;L^2(0,\ell))}
\leq \frac{\exp (T)r_0}{2\kappa_0}\,C_0^2 \, \Vert Q' \Vert^2_{L^2(0,T)},  \\ [8pt]
\displaystyle \Vert \varphi_{xx\tau} \Vert^2_{L^2(0,T;L^2(0,\ell))}
\leq \frac{\left (\exp (T)-1\right )r_0}{2\kappa_0}\,C_0^2 \, \Vert Q' \Vert^2_{L^2(0,T)}.
\end{array} \right.
\end{eqnarray*}
Further, as a consequence of (\ref{41}), we deduce that
\begin{eqnarray*}
\left. \begin{array}{ll}
\displaystyle \Vert Q' \Vert^2_{L^2(0,T)} \le \frac{4C_1^2}{\kappa_0} \Vert F\Vert^2_{L^2(\Omega_T)}+ 2\left [\Vert \theta'_0\Vert^2_{L^2(0,T)}+\Vert \theta'_\ell \Vert^2_{L^2(0,T)}\right ].
\end{array} \right.
\end{eqnarray*}
We obtain the necessary estimates (\ref{41}) by substituting this in the aforementioned inequalities. \hfill$\Box$

We can now use Theorem \ref{Theorem-3} and the estimates (\ref{11}) to justify the formal gradient formula (\ref{30}).

\begin{thm}\label{Theorem-4}
Suppose that conditions of Theorem \ref{Theorem-3} are satisfied. Then the Tikhonov functional  introduced in (\ref{20}) is Fr\'{e}chet differentiable, and for the Fr\'{e}chet gradient of this functional, the gradient formula (\ref{30}) is valid.
\end{thm}
{\bf Proof.} Denote by $\delta u(x,t)$ the solution of problem (\ref{1}) with $F(x,t)$ replaced by $\delta F(x,t)$. Then, it follows from the estimates (\ref{11}) applied to the solution $\delta u(x,t)$ that
\begin{eqnarray*}
\left. \begin{array}{ll}
\displaystyle \Vert \delta u_{x}(0,\,\cdot) \Vert^2_{L^2(0,T)}
\leq  \frac{C_1^2}{r_0}\, \Vert \delta F \Vert^2_{L^2(\Omega_T)},  \\ [10pt]
\displaystyle \Vert \delta u_{x}(\ell,\,\cdot) \Vert^2_{L^2(0,T)}
\leq  \frac{C_1^2}{r_0}\, \Vert \delta F \Vert^2_{L^2(\Omega_T)}.
\end{array} \right.
\end{eqnarray*}
Hence, the last two right-hand side integrals in (\ref{29}) are of the order
 $\mathcal{O}\left ( \Vert \delta F \Vert_{L^2(\Omega_T)}^2\right )$, which implies that
 the Tikhonov functional is Fr\'{e}chet differentiable and the formula (\ref{30}) is well-defined.
 This implies the proof. \hfill$\Box$

\section{The Lipschitz continuity of the Fr\'{e}chet gradient and monotonicity of iterations}

 It is a well-known fact that important properties, such as the monotonicity of iterations and the rate of convergence in gradient methods, are results of the Lipschitz continuity of the Fr\'{e}chet gradient $J'(F)$, $F \in \mathcal{F}$, of the Tikhonov functional \cite{Hasanov:2007}. Thus, for example, in the Landweber iteration algorithm $F^{n+1}=F^{n}-\omega_n J'(F^{n})$, applied to the inverse problem (\ref{1})-(\ref{2}), the relaxation parameter $\omega_n>0$ can be estimated through the Lipschitz constant \cite{Hasanov:Romanov:2021, Hasanov:2007}.

\begin{thm}\label{Theorem-5}
Assume that conditions of Theorem \ref{Theorem-3} are satisfied. Then the Fr\'{e}chet gradient of the Tikhonov functional, defined (\ref{30}), is Lipschitz continuous,
\begin{eqnarray}\label{44}
\Vert J'(F+\delta F)- J'(F)\Vert_{L^2(\Omega_T)}
\leq  L_G \,\Vert \delta F \Vert_{L^2(\Omega_T)},
\end{eqnarray}
with the Lipschitz constant
\begin{eqnarray}\label{45}
  \displaystyle L_G = \sqrt {\frac{\exp (T)-1}{2\kappa_0}} \,\ell^2 \,C_0 \,C_1,
\end{eqnarray}
where the constants $C_1>0$ and $C_0>0$ are defined in (\ref{12}) and (\ref{34}), respectively.
\end{thm}
{\bf Proof.} By the definition,
\begin{eqnarray}\label{46}
\Vert J'(F+\delta F)- J'(F)\Vert^2_{L^2(\Omega_T)}= \Vert \delta \varphi \Vert^2_{L^2(0,T;L^2(0,\ell))},
\end{eqnarray}
where $\delta \varphi (x,t):= \varphi (x,t;F+\delta F)-\varphi (x,t;F)$. Here, $\varphi (x,t;F+\delta F)$ and $\varphi (x,t;F)$ are the weak solutions of the adjoint problem (\ref{27}) for given $F +\delta F \in \mathcal{F}$ and $F \in \mathcal{F}$, respectively, and $\delta \varphi (x,t)$ solves the following problem:
\begin{eqnarray}\label{47}
\left\{ \begin{array}{ll}
\rho_A(x) \delta \varphi_{tt}-\mu(x)\delta \varphi_{t}-(T(x)\delta \varphi_{x})_{x}+ (r(x)\delta \varphi_{xx}-\kappa(x)\delta \varphi_{xxt})_{xx} =0,\\ [2pt]
\qquad \qquad \qquad \qquad \qquad \qquad \qquad \qquad (x,t)\in \Omega_{T}:=(0,\ell)\times (0,T);\\ [4pt]
\delta \varphi(x,T)=\delta \varphi_{t}(x,T)=0, \, x \in (0,\ell); \\ [4pt]
\delta \varphi(0,t)=0,\,
\left(r(x)\delta \varphi_{xx}-\kappa(x)\varphi_{xxt}\right)_{x=0}=\delta u_x(0,t); \\ [2pt]
\quad  \varphi(\ell,t)=0,\,\left (r(x)\varphi_{xx}-\kappa(x)\varphi_{xxt}\right)_{x=\ell}=
\delta u_x(\ell,t;F),\,t \in [0,T],
\end{array} \right.
\end{eqnarray}
with the inputs $\delta u_x(0,t):=u_x(0,t;F+\delta F)-u_x(0,t;F)$ and
$\delta u_x(\ell,t):=u_x(\ell,t;F+\delta F)-u_x(\ell,t;F)$. From the second estimate
of (\ref{33}) applied to the solution $\delta \varphi (x,t)$ of the problem (\ref{47})
we deduce that
\begin{eqnarray}\label{48}
\displaystyle \Vert \delta \varphi_{xx} \Vert^2_{L^2(0,T;L^2(0,\ell))} \qquad \qquad \qquad \qquad \qquad \qquad \qquad \qquad \qquad \qquad  \nonumber \\ [1pt]
\leq \left (\exp (T)-1\right ) C_0^2 \,
\left [\Vert \delta u_{xt}(0,\cdot) \Vert^2_{L^2(0,T)}+\Vert \delta u_{xt}(0,\cdot) \Vert^2_{L^2(0,T)} \right ].
\end{eqnarray}
We employ the second and fourth trace estimates in (\ref{11}) for the weak solution $\delta u(x,t):=u(x,t;F+\delta F) -u(x,t;F)$ of the problem (\ref{1}) with $F(x,t)$ replaced by $\delta F(x,t)$. Then, we have
\begin{eqnarray*}
\left. \begin{array}{ll}
\displaystyle \Vert \delta u_{xt}(0,\,\cdot) \Vert^2_{L^2(0,T)}
\leq  \frac{C_1^2}{\kappa_0}\, \Vert \delta F \Vert^2_{L^2(\Omega_T)},\\ [10pt]
\displaystyle \Vert \delta u_{xt}(\ell,\,\cdot) \Vert^2_{L^2(0,T)}
\leq  \frac{C_1^2}{\kappa_0}\,\Vert \delta F \Vert^2_{L^2(\Omega_T)}.
\end{array} \right.
\end{eqnarray*}
Hence, by (\ref{48}),
\begin{eqnarray}\label{49}
\displaystyle \Vert \delta \phi_{xx} \Vert^2_{L^2(0,T;L^2(0,\ell))}
\leq \frac{2\left (\exp (T)-1\right )}{\kappa_0} \,C_0^2 C_1^2 \,
\Vert \delta F \Vert^2_{L^2(\Omega_T)}.
\end{eqnarray}

On the other hand, by the inequality (\ref{13}), applied to the solution of (\ref{47}),
\begin{eqnarray*}
\displaystyle \Vert \delta \varphi_{x} \Vert^2_{L^2(0,T;L^2(0,\ell))}
\leq  \frac{\ell^2}{2}\,  \Vert \delta \varphi_{xx} \Vert^2_{L^2(0,T;L^2(0,\ell))}.
\end{eqnarray*}
If we consider here the inequality
\begin{eqnarray*}
\displaystyle \Vert \delta \varphi \Vert^2_{L^2(0,T;L^2(0,\ell))}
\leq  \frac{\ell^2}{2}\,  \Vert \delta \varphi_{x} \Vert^2_{L^2(0,T;L^2(0,\ell))},
\end{eqnarray*}
which can be easily proven for the function $\delta \varphi (x,t)$ satisfying the condition $\delta \varphi (0,t)=0$, then we obtain
\begin{eqnarray*}
\displaystyle \Vert \delta \varphi \Vert^2_{L^2(0,T;L^2(0,\ell))}
\leq  \frac{\ell^4}{4}\,  \Vert \delta \varphi_{xx} \Vert^2_{L^2(0,T;L^2(0,\ell))}.
\end{eqnarray*}
Taking into account this inequality in the estimate (\ref{49}), we arrive at the requested result.
\hfill$\Box$

%



\end{document}